\documentclass[12pt]{article}%
\usepackage{amssymb}
\usepackage{amsfonts}
\usepackage{graphicx}
\usepackage{amsmath}
\usepackage{appendix}
\usepackage{setspace}  
\providecommand{\U}[1]{\protect \rule{.1in}{.1in}}
\newenvironment{proof}[1][Proof. ]{\textbf{#1} }{\  \rule{0.5em}{0.5em}}
\newtheorem{Propo}{\textbf {Proposition}}
\newtheorem{lemma}{\textbf {Lemma}}
\newtheorem{Theorem}{\textbf {Theorem}}
\newtheorem{Assumption}{\textbf {Assumption}}
\newdimen \dummy
\dummy=\oddsidemargin
\begin{document}

\title{Remaining loads in a PH/M/$c$ queue with impatient customers\\
\bigskip
{\normalsize
Anders Rygh Swensen,\\Department of Mathematics, University of Oslo, \\P.O.
Box 1053, N-0316 Blindern, Oslo, Norway. \\Email:\ swensen@math.uio.no}}
\author{
\thanks{Department of Mathematics, University of Oslo, P.O.
Box 1053, N-0316 Blindern, Oslo, Norway. Email:\ swensen@math.uio.no} }
\date{}
\maketitle
\vspace{5mm}%


\vspace*{2cm}

\setcounter{page}{1}

\begin{abstract}
There are two contributions to the literature in this note. Firstly,
we point out that under some weak conditions the result in Swensen (1986) 
on the remaining loads in a GI/M/$c$ queue with impatient customers, 
derived under the assumption that the distribution of  inter-arrival times is
Coxian,  is also valid for  the more general phase type  distribution. 
From this result the distributions of the virtual and the actual waiting time can
easily be obtained. Secondly, the relation to an alternative expression for the distribution
of the virtual waiting time  derived by Kawanishi and Takine (2015) is also discussed.
\end{abstract}

\newpage

\section{Introduction}\label{intro} 
A queue where the customers that cannot be served within
a fixed time leave the system is usually referred to as a queue with 
impatient customers or a queue with
bounded waiting time and is often denoted as GI/M/c+D when the
arrival process is a renewal process.
To describe such systems balking and
reneging have also been used, see Ancker and Gafarian (1963). The motivation for
studying such queues has often been problems in telecommunication. The emergence of call
centers is an additional reason for the  attention  queues with impatient customers
has received.

We will consider a first-come first-served queue with $c$ servers where the service times are independent exponentially distributed random variables
with expectation $1/\mu$, the arrivals of new customers  follow a  renewal arrival process where
the distribution of the inter-arrival times  is of the phase type  and customers whose
waiting time exceeds the fixed impatience time $\tau$  leave the system.

One  message of the present note is that in a  PH/M/$c$  queue with impatient 
customers  explicit solutions, up to some matrix inversions, exist for the stationary 
asymptotic distribution of the 
remaining loads and consequently for  virtual  and the actual waiting time. It turns out that the procedure
from Swensen (1986) with a Coxian distribution describing  the renewal input  
process,  can be extended, with some adjustments, to the
more general case where a phase type distribution is used. The crucial result for the
generalization can be found in Lemma \ref{lem:lem1} in Appendix B.

Relevant papers dealing with impatient customers are 
Choi, Kim and Zhu (2004) and Kawanishi and Takine (2016) who  
considered the more general Markovian arrival process,  MAP/M/$c$ queue with bounded waiting 
or impatience time. He and Wu (2020) treated the situation where the
impatience time is not constant.
Kawanishi and Takine (2015) derived the
stationary distribution for the virtual waiting time in a
PH/M/$c$  queues with impatient customers.
It is worth noticing  that the results are complementary to the results in the
present paper.
Whereas the distributions of the waiting times, actual and virtual, in Swensen (1986) are
derived as an implication of the distribution of the loads of the servers, in 
Choi, Kim and Zhu (2004),  Kawanishi and Takine (2016) and Kawanishi and Takine (2015) the distribution of the virtual waiting time is derived directly by solving a 
second order differential equation. Thus, the results and methods of proof are different and a natural question is how they are related. We will discuss this issue. In particular we compare  the density for the virtual waiting time
derived by Kawanishi and Takine (2015) with the expression of the present note.

Although the distribution of the virtual waiting time is important, knowledge of the  simultaneous 
distribution of the loads of the servers represents an additional value. For example,
one can focus on the simultaneous distribution of a particular group of servers, cf. the famous Palm-Jacob{\ae}us  formula in telecommunications.
It is also worth remarking that the methods for the numerical computations of the distributions in the two cases are quite different. In  papers 
by Choi, Kim and Zhu (2004),  Kawanishi and Takine (2016) and Kawanishi and Takine (2015) 
the expressions that are used involve matrix-exponentials. On the other hand, 
the formula for the probability density of the remaining loads in
the present paper is quite easy to evaluate when a procedure for
finding matrix eigenvalues and eigenvectors is available.   

A phase type distribution can be described by the pair 
$(\boldsymbol{\gamma},T)$  where 
$\boldsymbol{\gamma}=
(\gamma_{1},\ldots,\gamma_{m})$  is
a probability vector and $T$ is
a non-singular $m \times m$ matrix with  negative diagonal elements and 
nonnegative off-diagonal elements.
The Coxian distributions define   the subclass of the phase type distributions where
$\boldsymbol{\gamma} =(1,0,\ldots,0)$ and the
only nonnegative elements of the matrix $T=\{t_{ij}\}$ are the elements
$t_{i(i+1)},\; i=1,\ldots,m-1$.

The note is organized as follows.
The distribution of the remaining loads  for  phase type distributed arrival times is derived in the next section.  In section 3 the relation between the implied distribution of the virtual
waiting time relying on the distribution of the remaining load and the more directly derived version in   
Kawanishi and Takine (2015)
is discussed.  In section 4   the solution of a linear equation which is important
for finding the asymptotic distribution of the remaining loads is considered.
In the appendices the necessary modifications of the Kolmogorov forward 
equations for the remaining loads for phase type distributed inter-arrival times 
are explained. The 
stationary solution to  these equations is found and
some lemmas crucial for the derivation are proved.
 
The arguments from the treatment of the
Coxian case carry over to a large extent  to the phase type case. Therefore there will be some
overlap with  Swensen (1986) when the arguments are verified in the 
phase type setup.
 
In the following we consider the scalar field of complex numbers and matrices
whose elements also are complex numbers.

\section{Distribution of remaining loads}\label{mainresult} 
Consider the process $(Y(t),V_1(t),\ldots,V_c(t))$ where $V_i(t)$ 
is the unfinished work, i.e. the remaining load,  for sever $i\;=1,\dots,c$ and
$Y(t)$ is the phase of the arrival process. Under the assumptions that the service times are
independent, exponentially distributed and that the arrival process is a renewal process
where the distribution of the inter-arrival times is of the
phase type, the process is Markovian.

The phase type distribution of the inter-arrival times satisfies
\begin{Assumption}
\label{ass:assump1}
\begin{itemize}
\item[i)]The m solutions $\{\eta_1,\ldots,\eta_m\}$ of 
$\det [c\mu \mathbf{e}\cdot \boldsymbol{\gamma} +T-c\eta I_m] = 0$
are distinct and are not equal to any of the eigenvalues of $T/c$.
\item[ii)]The matrix $T+\mathbf{T}\cdot \boldsymbol{\gamma} $
where $\mathbf{T}=-T\mathbf{e}, \mathbf{e}=(1,\cdots,1)'$, 
is irreducible and has
distinct eigenvalues.
\end{itemize}
\end{Assumption}

{\em Remark 1.} The claim in Assumption \ref{ass:assump1}
that there are no equal eigenvalues of $\mu \mathbf{e}\cdot \boldsymbol{\gamma} +T/c$
and no equal eigenvalues of $T+\mathbf{T}\cdot \boldsymbol{\gamma} $
simplifies the derivations since there is only one eigenvector associated
with each eigenvalue. We therefore  exclude the situation where the dimension of
the eigenspace of an eigenvalue is strictly less than the algebraic multiplicity 
of the eigenvalue, i.e. where the geometric multiplicity is less than the 
algebraic multiplicity, see Horn and Johnsen (2013).\\

The $m$ solutions $\{\eta_1,\ldots,\eta_m\}$
are  sorted in decreasing order according to size of the moduli. For a
conjugate pair the number with positive imaginary part is counted
first. Also, define
B as the matrix where the rows
are the left normalized eigenvectors of $T+ \mathbf{T}\cdot \boldsymbol{\gamma}$.
Let
D  be the diagonal matrix with diagonal elements
$(c\mu-\kappa_{\ell})^{-1},\;\ell=1,\ldots,m$ where 
$\kappa_{\ell},\;\ell=1,\ldots,m$ are the eigenvalues of 
$T+ \mathbf{T}\cdot \boldsymbol{\gamma}$
and where one, $\kappa_1$ say, is equal to $0$.
Let E  be the diagonal matrix with diagonal elements
$\exp(-\tau c \eta_{\ell}),\;\ell=1,\ldots,m$.

The following theorem is a description of the density of the
asymptotic distribution of the remaining loads.

\begin{Theorem}
\label{th:thm1}
Consider a first-come first-served GI/M/c queuing model
with waiting time bounded by $\tau$, service rate $\mu$
and inter-arrival distribution of the phase type $(\boldsymbol{\gamma},T)$.
Under Assumption \ref{ass:assump1}  the stationary 
asymptotic  distribution of 
the remaining loads $V_1,\ldots, V_c$ has density
\begin{equation}
\mu^{c-1}\exp(-\mu(v_1+\cdots+v_c) + c\mu (\tau \wedge \min v_i))
\sum_{k=1}^m \delta_k \underline{y}_{c}^k\mathbf{e}
\exp(-c\eta_{k}(\tau \wedge \min v_i))
\label{eq:eqtm1}
\end{equation}
for $\min V_i >0$
where $\boldsymbol{\delta}=(\delta_1,\ldots,\delta_m)$ is the solution
of 
\begin{equation}
\label{eq:eqtm2}
\boldsymbol{\delta}EY_{c-1}[c\mu\mathbf{T}\cdot \boldsymbol{\gamma}
B^{-1}DB-(c-1)\mu I_m+(c-1)\mu \mathbf{e}
\cdot\boldsymbol{\gamma} + T]=0
\end{equation}
and
\begin{equation}
\label{eq:eqtm3}
\sum_{k=1}^m  \delta_k \{ \sum_{i=0}^{c-1}\binom{c}{i}
\underline{y}_{i}^k \mathbf{e} + \underline{y}_{c}^k\mathbf{e}
[(\frac{1}{\mu}-\frac{1}{\eta_k})\exp{(-c\eta_k \tau) + \frac{1}{\eta_k}}]\}=1.
\end{equation}
The matrix $Y_{c-1}$  has rows
$\underline{y}_{c-1}^1,\ldots, \underline{y}_{c-1}^m$
where
$\underline{y}_{c-1}^k$  is the solution of the homogenous linear equation
\begin{align}
&\underline{y}_{c-1}^k[(c-1)\mu I_m -(c-1) \mu 
\mathbf{e}\cdot \boldsymbol{\gamma}-T- c \mu R_{k}] =0 \nonumber 
\end{align}
with $y_{c-1 1}^k =1$
and the $m \times m$ matrix $R_{k}
=\mathbf{T} \cdot \boldsymbol{\gamma} [c\eta_k I_m  -T]^{-1},\;
k=1,\ldots,m$.  

Furthermore, when $c=1$,  $\underline{y}_{1}^k=c\mu\underline{x}_{0}^kR_{k}$. 
When $c > 1$, 
$\underline{y}_{c}^k=c\mu \underline{y}_{c-1}^k R_{k}$
and
$\underline{y}_{0}^k,\ldots,\underline{y}_{c-2}^k,$ 
$k=1,\ldots,m$ are defined by the recursion
\begin{align}
&\underline{y}_{i}^k=(c-i)\mu \underline{y}_{i+1}^k( i\mu I_m - i \mu\mathbf{e}\cdot \boldsymbol{\gamma}-T)^{-1},
\; i=c-2,\ldots,1  \nonumber \\
&\underline{y}_0^k=-c\mu \underline{y}_1^k T^{-1}. \nonumber 
\end{align}
\end{Theorem}

{\em Remark 2.} From Lemma \ref{lem:lem2} it follows that the matrices $ [(c-1)\mu I_m -(c-1) \mu \mathbf{e}\cdot \boldsymbol{\gamma}\\-T- c \mu R_{k}]$ have rank m-1 
such that  $\underline{y}_{c-1}^k,\; k=1,\ldots,m$
are determined up to a normalizing constant.
Also it follows frem Lemma 1 in Neuts (1982) that  when $c>1$, the matrices
$( i\mu I_m - i \mu\mathbf{e}\cdot \boldsymbol{\gamma}-T),
\; i=1,\ldots,c-2$ 
are nonsingular. \\

The proof of the theorem is an elaboration of Swensen (1986)
for the case where the inter-arrival distribution is Coxian.
 In the first step to prove the main result  we show,
in Appendix A, that when $t \rightarrow \infty$ the functions
\begin{eqnarray}
q_{k0j} &=&y_{0j}^k,\;j=1,\ldots,m,     \label{eq:sol1}\\
q_{kij}(v_1,\ldots,v_i)&=&\mu^{i}y_{ij}^k\exp(-\mu(v_1+\cdots+v_i)),\;i=1,\ldots,
c-1,\;j=1,\ldots,m,  \label{eq:sol2}\\
q_{kcj}(v_1,\ldots,v_c)&=&\mu^{c-1}y_{cj}^k\exp(-\mu(v_1+\cdots+v_c)
+c(\mu-\eta_k)v_{(1)}),\;j=1,\ldots,m,  \label{eq:sol3}  
\end{eqnarray}
for $k=1,\ldots,m$,  
satisfy the differential equation defined by the Kolmogorov forward equations 
in the region $v_{(1)}=\min v_i \leq  \tau$. Since
$\eta_1,\ldots,\eta_m$ are distinct the solutions are independent
and the general solution in this
region therefore has the form $\mathbf{q} =\sum_{k=1}^m \delta_k \mathbf{q}_k $.

As a second step,
the function $\mathbf{q}$ is inserted in 
the Kolmogorov forward equation to
obtain a differential equation for values in the region $\min v_i > \tau$.
In appendix A it is shown that the solution has the form
\begin{equation}
K\exp(-\mu(v_1+\cdots+v_c))\boldsymbol{\gamma} B^{-1}D B
\nonumber
\end{equation}
where $K= \sum_{k=1}^m \delta_k c \mu^c\underline{y}_{c-1}^k
\mathbf{T}\exp(c(\mu-\eta_k)\tau)$.

Finally, in the third step the distribution of the remaining loads is determined by
requiring  that the  solutions  for $\min v_i \leq\tau$  and $\min v_i > \tau$
are  continuous at $\min v_i = \tau$. 

{\em Remark 3.}
The distribution of the virtual  and actual waiting time was derived 
for the $PH/M/c$ with impatient customers by 
Choi, Kim and Zhu (2004) and Kawanishi and Takine (2016).
But for other quantities the distribution of the remaining loads
is useful. An example is the distribution of the remaining
loads at $n$ particular servers. In Lemma A.7. in Swensen (1986) this is 
done  for  Coxian distributed inter-arrival times when 
$n=1$. Such a result can be generalized 
to the situation where $n>1$ and the arrival times 
are phase type distributed by a straight forward
application of the  the arguments of the 
present paper.

\section{Distribution of the virtual waiting time}
Building on the results in Choi et al. (2004)
Kawanishi and Takine (2015) found an explicit 
expression for the stationary  joint distribution of the phase of the
arrival process and the virtual waiting time.  
The density of the stationary virtual waiting time has the form
\begin{equation}
\mathbf{v}(v)\mathbf{e}=(v_1(v),\ldots,v_m(v))\mathbf{e}=p \mathbf{\hat{v}}\exp[(c\mu \mathbf{e}\cdot
\boldsymbol{\gamma} +T)(\tau-v) ]\mathbf{e},\; 0<v< \tau \label{eq:density2}\l
\end{equation}
where $\mathbf{\hat{v}}=(\hat{v}_1,\ldots,\hat{v}_m)= 
\boldsymbol{\gamma}(c\mu I-T)^{-1}/\boldsymbol{\gamma}(c\mu I-T)^{-1}\mathbf{e}$,
and $p$ is the normalizing constant
given in Kawanishi and Takine (2015).

However, an expression for the density of the stationary virtual
waiting time can also be obtained from the density of 
the remaining loads.
Differenciating the expression  for the complimentary
distribution function  in Swensen (1986) one gets the probability density of 
virtual waiting time when Assumption \ref{ass:assump1} is satisfied
\begin{equation}
f(v)= \sum_{k=1}^m\delta_k (y_c^k \mathbf{e})\;c\exp{(-c\eta_kv)}, \; 0<v <\tau \label{eq:density1}
\end{equation}
where  $\boldsymbol{\delta}=(\delta_1,\ldots,\delta_m)$ is defined in Theorem \ref{th:thm1}. By the same arguments as in the next section it follows 
that the expression (\ref{eq:density1}) has an imaginary part equal to zero.

The following proposition summarizes the relation
between the two expressions of the density.

\begin{Propo} If Assumption \ref{ass:assump1} holds the coefficients of
(\ref{eq:density1}) are given by
\begin{equation}
\delta_k (y_c^k \mathbf{e})= p\exp(c  \eta_k \tau )\sum_{j,\ell=1}^m \hat{v}_j f_{j,k} f^{k\ell}/c,\;k=1,\ldots,m \label{eq:prop1_01}
\end{equation}
where $c\eta_1,\ldots, c\eta_m $ are the eigenvalues of
$c\mu\mathbf{e}\cdot\mathbf{\gamma}+T$,
$F=\{f_{j\ell}\}_{j,\ell=1}^m$ 
is the matrix where the columns are the right eigenvectors of
the matrix $c\mu \mathbf{e}\cdot\mathbf{\gamma}+T$
and $F^{-1} = \{ f^{j\ell}\}_{j,\ell=1}^m$.

Also 
\begin{equation}
p=c\sum_{k=1}^m\delta_k (y_c^k\mathbf{e}) \exp(- c \eta_k \tau).
\label{eq:prop1_02}
\end{equation}
\end{Propo}

{\em Proof}. Assumption \ref{ass:assump1} $i)$ implies that 
$c\mu \mathbf{e}\cdot \boldsymbol{\gamma} +T$ is diagonalizable,   that is
$c\mu \mathbf{e}\cdot \boldsymbol{\gamma} +T= F \mbox{diag}(c\eta_1,\ldots,c\eta_m ) F^{-1}$
where the m eigenvalues $\{c\eta_1,\ldots,c\eta_m\}$
are  sorted in decreasing order according to the size of the moduli. 
From elementary properties of the
exponential matrix it  follows that
\begin{eqnarray*}
&&p \mathbf{\hat{v}}\exp[(c\mu \mathbf{e}\cdot
\boldsymbol{\gamma} +T)(\tau-v) ]\mathbf{e}  \\
&=&p \mathbf{\hat{v}}F\exp[ \mbox{diag} ({c\eta_1(\tau-v)},\ldots,{c\eta_m(\tau-v))}]F^{-1}\mathbf{e}\\
&=&p \mathbf{\hat{v}}F\mbox{diag} (\exp(c\eta_1(\tau-v)),\ldots,
\exp(c\eta_m(\tau-v)))F^{-1}\mathbf{e}\\
&=&p \sum_{k=1}^m[ \exp(c \tau \eta_k)\sum_{j,\ell=1}^m \hat{v}_j f_{j,k} f^{k\ell}]
\exp(-c\eta_k v),\; 0<v<\tau.
 \end{eqnarray*} 
Hence by comparing the coefficients
of $\exp(-c\eta_1),\ldots,  \exp(-c\eta_m), $
\[
\delta_k(y_c^k \mathbf{e})= p\exp( c\tau \eta_k )\sum_{j,\ell=1}^m \hat{v}_j f_{j,k} f^{k\ell}/c,\;k=1,\ldots,m
\]
which is (\ref{eq:prop1_01}). This shows       how the density (\ref{eq:density1}) for the virtual waiting time can be expressed using the
results (\ref{eq:density2})     from Kawanishi and Takine (2015). 

Furthermore, summing over $k$, since $\sum_{k=1}^m f_{j,k} f^{k\ell} =1$ if 
$j=\ell$ and $  0$ if $j\neq \ell$,
\begin{eqnarray*} 
\sum_{k=1}^m\delta_k (y_c^k \mathbf{e})\exp(- c\tau \eta_k )&=&
p\sum_{j,\ell=1}^m \hat{v}_j \sum_{k=1}^m f_{j,k} f^{k\ell}/c=
p\sum_{j=1}^m \hat{v}_j/c= p 
\mathbf{\hat{v}}\mathbf{e}/c
\\&=&
\frac{1}{c}
p\boldsymbol{\gamma}(c\mu I-T)^{-1}\mathbf{e}/\boldsymbol{\gamma}(c\mu I-T)^{-1}\mathbf{e}= p/c
\end{eqnarray*} 
which is (\ref{eq:prop1_02}).
\hspace{ 14 cm}$\blacksquare$

{\em Remark  4.}  For numerical computations the 
representation in (\ref{eq:density1}) yields an alternative
to the approach using matrix exponentials suggested by 
 Kawanishi and Takine (2015). With a procedure 
 computing eigenvalues and eigenvectors the  
 quantities $\delta_1,\ldots,\delta_m$ can be found  as
 described in Theorem \ref{th:thm1}.

\section{Solution of the linear equations (\ref{eq:eqtm2}) and
 (\ref{eq:eqtm3})} 
 
 The linear equations (\ref{eq:eqtm2}) and
 (\ref{eq:eqtm3}) can be solved directly by
 writing them on the form $A\boldsymbol{\delta}'=e_{m+1}$ where $A$ is
 a $(m+1)\times m$ matrix and $e_{m+1}$ is a $(m+1)$-dimensional
 vector with  all elements  equal to $0$ except the last which equals
 $1$. Then, using a $QR$ factorization, as described in Theorem 2.1.14
 in Horn and Johnson (2013), there is a   $(m+1)\times m$ matrix $Q$ 
 with orthonormal columns and an upper triangular 
 $m\times m$ matrix $R$ such that the linear equations may be 
 expressed as $QR\boldsymbol{\delta}'=e_{m+1}$ from which
 $\boldsymbol{\delta}$ is easily found.
 
 But there is an alternative way to obtain a solution starting with solving
 \begin{equation}
\label{eq:met2_1}
\boldsymbol{\phi}[c\mu\mathbf{T}\cdot \boldsymbol{\gamma}
B^{-1}DB-(c-1)\mu I_m+(c-1)\mu \mathbf{e}
\cdot\boldsymbol{\gamma} + T]=0
\end{equation}
which is determined up to a normalization  by
Lemma \ref{lem:lem3}.  Normalizing the first element of
$\boldsymbol{\phi}$ as $1$, it follows from
Lemma \ref{lem:lem0} that all elements of
$\boldsymbol{\phi}$ are real. 
\begin{lemma}
\label{lem:lem0}
The solutions of equation (\ref{eq:met2_1}) have
real elements.
\end{lemma}
{\em Proof.} If $\boldsymbol{\kappa}$ is the diagonal
matrix with diagonal elements $\kappa_1,\ldots,\kappa_m$,\\
$B(T+ \mathbf{T}\cdot \boldsymbol{\gamma})B^{-1}=
\boldsymbol{\kappa}$.
Furthermore, $B c\mu I_m B^{-1} = c\mu I_m$.
Subtracting $B(-T- \mathbf{T}\cdot \boldsymbol{\gamma}+c\mu I_m)B^{-1}=
c\mu I_m-\boldsymbol{\kappa}_m$. Then
$B(c\mu I_m -T- \mathbf{T}\cdot \boldsymbol{\gamma})^{-1}B^{-1}=
(c\mu I_m-\boldsymbol{\kappa}_m)^{-1}$
and $(c\mu I_m -T- \mathbf{T}\cdot \boldsymbol{\gamma})^{-1}=
B^{-1}(c\mu I_m-\boldsymbol{\kappa}_m)^{-1}B=B^{-1}DB$.
Thus since the matrix
$(c\mu I_m -T- \mathbf{T}\cdot \boldsymbol{\gamma})$ has real elements,
so has the matrix  $B^{-1}DB$ and the matrix
$c\mu\mathbf{T}\cdot \boldsymbol{\gamma}
B^{-1}DB-(c-1)\mu I_m+(c-1)\mu \mathbf{e}
\cdot\boldsymbol{\gamma} + T$. Therefore also
the solutions of (\ref{eq:met2_1}) must be real. $\blacksquare$

Remark that $\boldsymbol{\phi}=\boldsymbol{\delta}EY_{c-1}$.  
The matrix  $Y_{c-1}$  defined in Theorem 1 has rows corresponding to $\eta_1,\ldots,\eta_m$.
Then $\bar{Y}_{c-1}$  are the  rows corresponding to $\bar{\eta}_1,\ldots,\bar{\eta}_m$
where $\bar{\eta}_i$ is the conjugate of $\eta_i$ if this number is not
real.  Let $\boldsymbol{\delta}$ and $\boldsymbol{\tilde{\delta}}$
be the solutions of  equation (\ref{eq:eqtm2}). Since 
the elements  of $\boldsymbol{\phi}$ are real
\begin{equation}
\boldsymbol{\phi}=\boldsymbol{\delta}Y_{c-1}E= 
  \boldsymbol{\tilde{\delta}}\bar{Y}_{c-1}\bar{E}. \nonumber
\end{equation}
The matrices $\bar{Y}_{c-1}$ and $\bar{E}$ are nonsingular,
such that $\boldsymbol{\bar{\delta}}=  \boldsymbol{\tilde{\delta}} $.
This means that in the solutions of equation (\ref{eq:eqtm3})
the complex elements occur in conjugate pairs. But then 
\begin{eqnarray}
\sum_{k=1}^m  \delta_k     \{ \sum_{i=0}^{c-1}\binom{c}{i}
\underline{y}_{i}^k \mathbf{e} + \underline{y}_{c}^k\mathbf{e}
[(\frac{1}{\mu}-\frac{1}{\eta_k})\exp{(-c\eta_k \tau) + \frac{1}{\eta_k}}]\}\nonumber\\
=\sum_{k=1}^m  \delta_k     \underline{y}_{c-1}^k\{
\sum_{i=0}^{c-1}\binom{c}{i} 
\prod_{\ell=i}^{c-2} (c-\ell)\mu^{c-\ell+1} ( \ell \mu I_m -\ell \mu\mathbf{e} 
\cdot \boldsymbol{\gamma}-T)^{-1} \label{eq:compsum}
\\ +cI_m+
 c\mu R_{k}
[(\frac{1}{\mu}-\frac{1}{\eta_k})\exp{(-c\eta_k \tau) + \frac{1}{\eta_k}}]
\}\mathbf{e} \nonumber
 \end{eqnarray}
must be a real number, since conjugate pairs of $\delta_i$ correspond
to conjugate pairs of other part of the terms in the sum
(\ref{eq:compsum}). But the sum
of a complex number and its conjugate must be real. Therefore 
(\ref{eq:compsum}) is real and can be used for normalizing
$\boldsymbol{\delta}$. 

\section{ Appendices}
\appendix

\section{Stationary solutions of the Kolmogorov forward equations}
We will use the following notation  $T=\{ t_{jl }\}_{j,l=1}^m$, $t_{jj}=-\lambda_j$,
$t_{jl} = \lambda_j \beta_{jl},\: j\neq l $ for $j,l=1\ldots,m$.
Also $\beta_{jj}=0$ and $\alpha_j=1-\sum_{l=1}^m \beta_{jl} $
for $j=1,\ldots,m$ such that $\mathbf{T}=(\lambda_1\alpha_1,\ldots,
\lambda_m\alpha_m)'$.

The forward Kolmogorov equations are  a  generalization  of those 
used in Swensen (1986). Due to the restrictions on  $\boldsymbol{\gamma}$ and $T$ in the Coxian case
some  modifications are necessary when a
general phase type distribution is used to describe the arrival process.

Let $p_{ij}(t+h,v_1,\ldots,v_i)$ be the density of  the remaining load for
$i$ occupied servers, $i,\;i=0,\ldots,c$, and arrival process in state $j,\;j=1,\ldots,m$
at time $t+h$.
Letting first $t \rightarrow \infty$ and then $h \rightarrow 0$
it follows that the following equations
must be satisfied by the stationary solutions
\begin{eqnarray}
\label{eq:seceq1}
[c p_{11}(0),\ldots,c p_{1m}(0)]&=-&[\sum_{l=1}^m  p_{0l}t_{l1},
\ldots,
\sum_{l=1}^m  p_{0l}t_{lm}],
\end{eqnarray}
\begin{eqnarray}
\label{eq:seceq2}
-\sum_{l=1}^{i}\frac{dp_{ij}}{dv_l}
&=&
\sum_{l =1}^m
p_{il}(v_1,\ldots,v_i) t_{lj}\\
&+&\frac{\gamma_j}{(c-i+1)}
\sum_{l=1}^m \sum_{k=1}^{i}
p_{i-1l}(v_1,\ldots,v_{k-1}, v_{k+1},\ldots,v_i )
\lambda_l\alpha_{l}\mu \exp (-\mu v_k)
 \nonumber \\
 &+&(c-i) p_{i+1j}(v_1,\ldots,v_i,0)
 \nonumber
 \end{eqnarray}
for $i=1,\ldots,c-1$ when $ c>1 $
and
 \begin{eqnarray}
 \label{eq:seceq3}
 -\sum_{l=1}^{c}\frac{dp_{cj}}{dv_l}
&=&\sum_{l =1}^m 
p_{cl}(v_1,\ldots,v_c) t_{lj}  \\
&+&\gamma_j \sum_{l=1}^m \sum_{k=1}^{c}
p_{c-1l}(v_1,\ldots,v_{k-1}, v_{k+1},\ldots,v_c)
\lambda_l\alpha_{l} \mu \exp (-\mu v_k)
 \nonumber \\
 &+&\gamma_j \sum_{l=1}^m \sum_{k=1}^{c}\int_0^{\tau \wedge \min v_i}
p_{cl}(v_1,\ldots,v_{k-1}, v,v_{k+1},\ldots,v_c)
\lambda_l\alpha_{l}\mu \exp (-\mu (v_k-v))dv
 \nonumber \\
 &+&\gamma_j I(\tau < \min v_i) \sum_{l=1}^m p_{cl}(v_1,\ldots,v_{c})
 \lambda_l \alpha_l  \nonumber
 \end{eqnarray}
 for  $j=1,\ldots,m$. 
 
Firstly, we claim that a solution of the previous equations 
has the form defined by the equations
(\ref{eq:sol1})-(\ref{eq:sol3}) in the region $\min v_i \leq \tau$.
As in Swensen (1986) where the similar situation is treated for the Coxian distributions the verification is by straight forward insertion
in (\ref{eq:seceq1}) and (\ref{eq:seceq2}) for $i=1,\ldots,c-1$. 
 
For the case  $i=c$ the insertion of  
equations  (\ref{eq:sol2}) and (\ref{eq:sol3}) in 
the left- and right-hand side of 
equation
(\ref{eq:seceq3}) yields 
that the following equations must be satisfied 
for $k=1,\ldots,m$
\begin{equation}
\label{eq: ceq2}
\frac{\mu}{\mu-\eta_k} \underline{y}_{c}^k \mathbf{T} =
c\mu  \underline{y}_{c-1}^k \mathbf{T}
\end{equation}
and

\begin{equation} 
\underline{y}_c^k[ c\eta I_m -T ]= \frac{\mu [\underline{y}_c^k \mathbf{T}] \boldsymbol{\gamma} }{ \mu-\eta_k}. \label{eq: ceq1}
\end{equation}
It follows from Lemma \ref{lem:lem1} that $\eta_k$ is a solution of
$(\mu-\eta_k)/\mu= \boldsymbol{\gamma} 
[c\eta_k I_m -T]^{-1} \mathbf{T}
$.  Therefore the equation (\ref{eq: ceq2}) may be written
\begin{equation}
\label{eq: ceq3}
\underline{y}_c^k \mathbf{T} = c\mu \underline{y}_{c-1}^k \mathbf{T}\cdot \boldsymbol{\gamma}
[c \eta_k I_m - T]^{-1}\mathbf{T}= c \mu  \underline{y}_{c-1}^k R_{k} \mathbf{T},
\end{equation}
where the last equality follows from the definition of $R_{k}$. But
the definition of $\underline{y}_c^k$ and $\underline{y}_{c-1}^k$
implies (\ref{eq: ceq3}).

To verify (\ref{eq: ceq1}) note that from the definition of $\underline{y}_{c-1}^k$ and $R_{k}$ it follows that
\begin{equation*}
\underline{y}_c^k[c\eta_k I_m -T]= c\mu \underline{y}_{c-1}^k R_{k}[c\eta_k I_m -T]= c\mu  \underline{y}_{c-1}^k  \mathbf{T}\cdot 
\boldsymbol{\gamma}
\end{equation*}
such that 
the equations (\ref{eq: ceq1}) may be written, using the definition
of $\underline{y}_{c-1}^k$ and $\underline{y}_{c}^k$
\begin{equation}
\label{eq: ceq4}
c\mu \underline{y}_{c-1}^k \mathbf{T}\cdot \boldsymbol{\gamma}
= \frac{\mu}{(\mu-\eta_k)} \underline{y}_{c}^k\mathbf{T} \cdot \boldsymbol{\gamma}.
\end{equation}
But from the definition of $\underline{y}_{c}^k $,
$\underline{y}_{c-1}^k $   and Lemma \ref{lem:lem1} it follows that
\begin{equation}
\label{eq: aux1}
\underline{y}_c^k \mathbf{T} = c\mu \underline{y}_{c-1}^k \mathbf{T}\cdot \boldsymbol{\gamma}
[c \eta_k I- T]^{-1}\mathbf{T}= c\mu \underline{y}_{c-1}^k \mathbf{T}\frac{(\mu - \eta_k)}{\mu} 
\end{equation}
which implies  (\ref{eq: ceq4}) and therefore (\ref{eq: ceq1}).
Thus the proposed solutions given by equations (\ref{eq:sol1})-(\ref{eq:sol3})
satisfy  the requirements from equations (\ref{eq:seceq1})-(\ref{eq:seceq3})
in the area $\min v_i \leq \tau$ for each $\eta_1,\ldots,\eta_m$.

Secondly, the general solution will be a linear combination of these 
solutions. Inserting such a  solution in the equation
(\ref{eq:seceq3})  
one gets when $\min v_i > \tau$ the following equation on vector form 
\begin{equation}
\label{eq: diff01}
(-\sum_{l=1}^{c}\frac{dp_{c1}}{dv_l},\ldots,-\sum_{l=1}^{c}\frac{dp_{cm}}{dv_l})
-(p_{c1},\ldots,p_{cm})(T+\mathbf{T}\cdot \boldsymbol{\gamma})=K\boldsymbol{\gamma}\exp(-\mu(v_1+\ldots+v_c))
\end{equation}
where by using (\ref{eq: aux1})

\begin{equation*}
K=\sum_{k=1}^m \delta_k
\underline{y}_{c}^k\mathbf{T}\mu^c \exp(c(\mu-\eta_k)\tau))/(\mu - \eta_k)
= c\sum_{k=1}^m \delta_k \underline{y}_{c-1}^k\mathbf{T}\mu^c\exp(c(\mu-\eta_k)\tau)).
\end{equation*}

Let $\kappa_1,\ldots,\kappa_m$ be the eigenvalues of 
$ T + \mathbf{T}\cdot \boldsymbol{\gamma}$,
$B=(b_1',\ldots,b_m')'$ be  the matrix with  the  normalized left eigenvectors
of $ T + \mathbf{T}\cdot \boldsymbol{\gamma}$ as rows and let
D be the diagonal matrix with diagonal elements
$1/(c\mu-\kappa_1),\ldots,1/(c\mu-\kappa_m)$. Since 
$ T + \mathbf{T}\cdot \boldsymbol{\gamma}$ is a generator
the Perron-Frobenius theorem implies that one eigenvalue, $\kappa_1=0$ say,
is equal to zero and the rest have strictly negative real parts.
Also by Assumption 1 $ii)$ $\kappa_1,\ldots,\kappa_m$
are distinct. Since  $\kappa_2,\ldots,\kappa_m$ have negative real parts they are not
 equal to $c\mu > 0$.

As in Swensen (1986) one can show that
\begin{equation}
\nonumber
\psi_j(T + \mathbf{T}\cdot \boldsymbol{\gamma})=
b_j(T + \mathbf{T}\cdot \boldsymbol{\gamma})\exp(-\kappa_j v_1)\phi_j=
\kappa_j b_j \exp(-\kappa_j v_1)\phi_j =-\sum_{l=1}^c \frac{d\psi_j}{dv_l} 
\end{equation}
where $\phi_1,\ldots,\phi_m$ are arbitrary differentiable functions in one variable
and $\psi_j(v_1+\ldots+v_c)= b_j\exp(-\kappa_j v_1)\phi(-(c-1)v_1+v_2+\cdots+v_c),\;
j=1,\ldots,m$.
Therefore  $\psi_1,\ldots,\psi_m$  are $m$ independent
solutions 
of the homogenous differential equation corresponding to  (\ref{eq: diff01}).

Now define,
\begin{equation}
\nonumber
\psi_0(v_1,\ldots,v_c)= K\exp(-\mu(v_1+\cdots+v_c))\boldsymbol{\gamma}
B^{-1} DB.
\end{equation}
Then
\begin{equation}
\nonumber
-\sum_{l=1}^c \frac{d \psi_0}{d v_l}= c\mu K \exp(-\mu(v_1+\cdots+v_c))
\boldsymbol{\gamma}B^{-1} DB.
\end{equation}
Thus, if $\boldsymbol{\kappa}$ is the diagonal matrix with diagonal elements $\kappa_1,\ldots, \kappa_m$,
\begin{eqnarray}
\nonumber
-\sum_{l=1}^c \frac{d \psi_0}{d v_l} - \psi_0
(T + \mathbf{T}\cdot \boldsymbol{\gamma})
&=&K\exp(-\mu(v_1+\cdots+v_c))
\boldsymbol{\gamma} B^{-1} ( c\mu D -D \boldsymbol{\kappa})B\\
&=&
K\boldsymbol{\gamma}\exp(-\mu(v_1+\cdots+v_c))   \nonumber
\end{eqnarray}
such that $\psi_0$ is a particular solution of 
 (\ref{eq: diff01}) and the general solution is 
 $\psi_0 + \sum_{l=1}^m \nu_l \psi_l $.
 
 Since $\kappa_1=0$ and  $\kappa_2,\ldots,\kappa_m$
 have strictly negative real parts, $ \nu_1=\cdots=\nu_m=0$
 because $p_{cj} \rightarrow 0, \; j=1,\ldots m $ as $\sum_{l=1}^c v_l^2 \rightarrow \infty.$
 
 Finally, from the Kolmogorov forward equations it follows that 
 \begin{equation} 
 \nonumber
 p_{cj}(v_1,\ldots,v_c)= \lim_{h \rightarrow 0 } p_{cj}(v_1+h,\ldots,v_c+h), 
 \; j=1,\ldots m 
 \end{equation}
 such that for $(v_1,\ldots,v_c) \in \{v_1,\ldots,v_c) : \min v_i =\tau  \}$
 \begin{equation} 
 \nonumber
 \sum_{k=1}^m \delta_k c\mu \underline{y}_{c-1}^k\exp(c(\mu- \eta_k)\tau)
 \mathbf{T}\cdot
 \boldsymbol{\gamma}B^{{-1}}DB
= \sum_{k=1}^m \delta_k \underline{y}_{c}^k\exp(c(\mu- \eta_k)\tau).
 \end{equation} 
 Using the definition of $\underline{y}_{c}^k$ and letting $Y_{c-1}$
 be the matrix with rows $\underline{y}_{c-1}^k,\;k=1,\ldots,m$ 
 this equals 
 \begin{equation} 
 \label{eq:fisol}
 \boldsymbol{\delta } EY_{c-1}[c \mu \mathbf{T}\cdot
 \boldsymbol{\gamma} B^{-1}DB- (c-1)\mu I_m
 + (c-1) \mu \mathbf{e} \cdot \boldsymbol{\gamma} +T]=0.
  \end{equation}
 The matrix $E$  has full rank and so has $Y_{c-1}$
 by Lemma \ref{lem:lem4}. Hence it follows from 
 Lemma \ref{lem:lem3} that equation (\ref{eq:fisol}) has a solution 
 up to a multiplicative constant, which is found by normalization.
 
\section{Some lemmas}
\begin{lemma}
\label{lem:lem1}
Under Assumption 1 $i)$ the set of solutions 
$\eta_1,\ldots,\eta_m$
of the equation 
\begin{equation}
\det[c\mu \mathbf{e}\cdot \boldsymbol{\gamma} +T-c\eta I_m]=0 \nonumber
\end{equation}
is contained in the set of solutions of 
\begin{equation}
\frac{\mu -\eta}{\mu}=\boldsymbol{\gamma} [ c\eta I_m -T]^{-1} \mathbf{T}.
\nonumber
\end{equation}
\end{lemma}

{\em Proof.} 
First note that from
Assumption 1 i) $\eta$ is different from the 
eigenvalues of $T/c$ so that $(c\eta I_m -T)$ is
invertible. Then
\begin{eqnarray}
\nonumber
f(c\eta)&=& \boldsymbol{\gamma}(c\eta I_m -T)^{-1} (-T) \mathbf{e}=
\boldsymbol{\gamma}(c\eta I_m -T)^{-1} (c\eta I_m-T-c\eta I_m) \mathbf{e}\\
\nonumber
&=&\boldsymbol{\gamma}\mathbf{e} - c\eta \boldsymbol{\gamma}(c\eta I_m -T)^{-1}  \mathbf{e}=1-c\eta\boldsymbol{\gamma}(c\eta I_m -T)^{-1} \mathbf{e}.
\end{eqnarray}
Using Cauchy's formula for the determinant of a rank-one perturbation,
see Horn and Johnson (2013) p. 26,   the relation
\begin{equation}
\nonumber
0=\det[c\mu \mathbf{e}\cdot \boldsymbol{\gamma} +T-c\eta I_m]=
\det(T-c\eta I_m ) ( 1 + c\mu \boldsymbol{\gamma}(T-c\eta I_m )^{-1} \mathbf{e})
\end{equation}
implies that 
\begin{equation}
\nonumber
 \frac{1}{\mu}=- c\boldsymbol{\gamma}(T-c \eta I_m )^{-1} \mathbf{e}.
\end{equation}
Hence 
\begin{equation}
\frac{\mu -\eta}{\mu}=1-\frac{\eta}{\mu}= 1- c\eta \boldsymbol{\gamma}
(c \eta I_m  -T)^{-1} \mathbf{e} = f(c\eta)
=\boldsymbol{\gamma}(c \eta I_m -T)^{-1} \mathbf{T}. \;\;\;\nonumber \blacksquare
\end{equation}

\begin{lemma}
\label{lem:lem2}
The matrices 
\begin{equation}
(c-1)\mu I_m -(c-1) \mu \mathbf{e}\cdot \boldsymbol{\gamma}-T- c \mu R_{k},\;
k=1,\ldots,m
\label{eq:lem2eq} 
\end{equation}
have rank m-1.
\end{lemma}
{\em Proof.}
As pointed out in Swensen (1986) it suffices to prove that 
the matrix (\ref{eq:lem2eq}) is singular.
The case where $R_{k}-I$ is non-singular is proved as in Swensen (1986).
To deal with the case where $R_{k}-I$  is singular note that
from Lemma \ref{lem:lem1} it follows that 
$R_{k}\mathbf{T}  =\mathbf{T} \cdot \boldsymbol{\gamma}[c\eta_k I_m  -T]^{-1}\mathbf{T}
=[(\mu-\eta_k)/\mu] \mathbf{T}$ so $\mathbf{T}$ is a right eigenvector of $R_{k}$
and $(\mu-\eta_k)/\mu$ is the non-zero eigenvalue. When 
$R_{k}-I$ is singular, $1$  has to be the non-zero eigenvalue of
$R_{k}$ such that $1=(\mu-\eta_k)/\mu  $ and $\eta_k=0$. 
Furthermore, for $\eta_k=0$ and using
Cauchy's formula for the determinant of a rank-one perturbation
\begin{equation}
\nonumber
0=\det (c \mu \mathbf{e} \cdot \boldsymbol{\gamma} +T-c\eta_k I_m)=
\det (c \mu \mathbf{e} \cdot \boldsymbol{\gamma} +T)=
(\det{T} )(1+c\mu \boldsymbol{\gamma} T{^{-1}} \mathbf{e} )
\end{equation}
such that $c\mu \boldsymbol{\gamma} T{^{-1}} \mathbf{e}=-1$
and 
\begin{equation}
[(c-1)\mu I_m -(c-1) \mu \mathbf{e}\cdot \boldsymbol{\gamma}-T- c \mu R_{k}]\mathbf{e} =( -T+ c \mu \mathbf{T} \cdot \boldsymbol{\gamma} T^{-1})\mathbf{e}
=(-T\mathbf{e}) -\mathbf{T})=0, \nonumber
\end{equation}
which shows that the matrix  (\ref{eq:lem2eq})
must be singular also when $R_{k}-I$  is singular.
$\blacksquare$

\begin{lemma}
\label{lem:lem3}
Let $\kappa_1,\ldots,\kappa_m$ be the eigenvalues of 
$ T + \mathbf{T}\cdot \boldsymbol{\gamma}$ and
B the matrix with rows as normalized left eigenvectors
of $ T + \mathbf{T}\cdot \boldsymbol{\gamma}$ and let
D be the diagonal matrix with diagonal elements
$1/(c\mu-\kappa_1),\ldots,1/(c\mu-\kappa_m)$.
Under Assumption 1 $ii)$ the matrix 
\begin{equation}
c\mu\mathbf{T}\cdot \boldsymbol{\gamma} B^{-1}DB  -
(c-1)\mu I_m +(c-1) \mu \mathbf{e}\cdot  \boldsymbol{\gamma}     +T
\nonumber 
\end{equation}
has rank m-1.
\end{lemma}
{\em Proof.}
From Lemma 1 in Neuts (1982) it follows that
the matrix
$(c-1)\mu I_m -(c-1) \mu \mathbf{e}\cdot \boldsymbol{\gamma} -T$ has rank m.
The matrix
$c\mu\mathbf{T}\cdot \boldsymbol{\gamma} B^{-1}DB  -
(c-1)\mu I_m +(c-1) \mu \mathbf{e}\cdot \boldsymbol{\gamma} +T$ therefore has rank
at least m-1. Hence, it suffices to show that it is singular.

Let $\boldsymbol{\kappa}$ be the diagonal  matrix with
diagonal elements $\kappa_1,\ldots,\kappa_m$.
Without loss of generality  we can assume $\kappa_1=0$.
Then 
\begin{eqnarray}
B(T+ \mathbf{T}\cdot \boldsymbol{\gamma})\mathbf{e} = \boldsymbol{\kappa}B \mathbf{e},
\nonumber \\
BT\mathbf{e} + B\mathbf{T}\cdot \boldsymbol{\gamma}\mathbf{e} = \boldsymbol{\kappa}B \mathbf{e},  \nonumber \\
0=-B\mathbf{T} + B\mathbf{T} =\boldsymbol{\kappa}B \mathbf{e}.  \nonumber 
\end{eqnarray}
But $\kappa_i \neq 0$ such that $b_i\mathbf{e} =0$ for  the ith row of B, $b_i$,\;
$i=2,\ldots,m$ and
$B\mathbf{e}= (\mathbf{e}' b_1',0,\cdots,0)' =a(1,0,\cdots,0)'=a \mathbf{e}_1$.
This implies that $\mathbf{e}=B^{-1}B \mathbf{e}=aB^{-1}\mathbf{e}_1$ and
\begin{equation}
1=\boldsymbol{\gamma}\mathbf{e} =a\boldsymbol{\gamma}B^{-1}\mathbf{e}_1.
\label{eq:lem3_1}
\end{equation}
Also, 
\begin{equation*}
DB\mathbf{e}  = Da\mathbf{e} _1 = \frac{a}{c\mu}\mathbf{e}_1.
\label{eg:lem3_2}
\end{equation*}
Hence,
\begin{eqnarray}
c\mu\mathbf{T}\cdot \boldsymbol{\gamma} B^{-1}DB \mathbf{e}
&=&c\mu\mathbf{T}\cdot \boldsymbol{\gamma} B^{-1} \frac{a}{c\mu}\mathbf{e}_1
\nonumber\\
&=&c\mu\mathbf{T}\cdot (a\boldsymbol{\gamma} B^{-1} \mathbf{e}_1)\frac{1}{c\mu}
=\mathbf{T}= -T \mathbf{e}\nonumber
\end{eqnarray}
by (\ref{eq:lem3_1}), which shows that 
$c\mu\mathbf{T}\cdot \boldsymbol{\gamma} B^{-1}DB  -
(c-1)\mu I_m +(c-1) \mu \mathbf{e}\cdot \boldsymbol{\gamma} +T$
is singular.
$\blacksquare$
 
\begin{lemma}

\label{lem:lem4}
The vectors $y_{c-1}^k,\; k=1\ldots,m$ are linearly  independent.
\end{lemma}
{\em Proof.} See Lemma A.6 in Swensen (1986). The argument
for the Coxian distributions is also valid for phase type distributions.
$\blacksquare$


\begin{thebibliography}{99}
\footnotesize
\bibitem{ref8} C.J. Ancker,  A.V. Gafarian, Some queuing problems
with balking and reneging, Oper. Res. 11 (1963) 88-100.
\bibitem{ref1}
B. Choi, B. Kim, D. Zhu,  MAP/M/c Queue with constant impatient time, Math.  Oper. Res. 29 (2004) 309-325. 
\bibitem{he}
Q-M. He,  H. Wu, Multi-Layer MMFF Processes and the MAP/PH/K+GI Queue: Theory and Algorithms,
Queueing Models and Service Management 3 (2020) 37-87.
\bibitem{ref7} R.A. Horn, C.R. Johnson, Matrix Analysis, sec. ed.,
Cambridge University Press, New York,  (2013).  
\bibitem{ref2}
K. Kawanishi, T.  Takine,  A note on the virtual waiting time
in the stationary PH/M/c+D queue, J. Appl.  Probab. 52 (2015) 899-903.   
\bibitem{ref3}
K. Kawanishi, T. Takine, MAP/M/c and  M/PH/c queues
with constant impatience times,  Queuing Syst. 82 (2016) 381-420.
\bibitem{ref6}  M.F. Neuts, Explicit steady-state solutions to some elementary queuing models, Oper. Res. 30 (1982) 480-489.
\bibitem{ref4}  A.R. Swensen, On a GI/M/c queue with bounded waiting times,  Oper. Res. 34 (1986) 895-908.
\end{thebibliography}
\end{document}